\documentclass[12pt]{amsart}
\usepackage{amssymb}
 
\newtheorem{thm}{Theorem}[section]

\newtheorem{lem}{Lemma}[section]
\newtheorem{cor}{Corollary}[section]
\begin{document}
\author{G\'erard Endimioni}
\address{C.M.I-Universit\'{e} de Provence,
39, rue F. Joliot-Curie, F-13453 Marseille Cedex 13}
\email{endimion@cmi.univ-mrs.fr}
\title[Polynomial Functions]{On the Group of Polynomial Functions in a Group}
% \keywords{}
% \dedicatory{}
%\thanks{ }

\begin{abstract}
Let $G$ 
be a group and let $n$ be a positive integer.
A polynomial function in $G$ is a function from $G^n$ to $G$ of the form
$(t_{1},\ldots,t_{n})\to f(t_{1},\ldots,t_{n})$, where $f(x_{1},\ldots,x_{n})$
is an element of the free product of $G$ and the free group of rank $n$
freely generated by $x_{1},\ldots,x_{n}$. There is a natural 
definition for the product of two polynomial functions; equipped with this 
operation, the set  $\overline{G}[x_{1},\ldots,x_{n}]$ of polynomial 
functions is a group. We prove that this group is polycyclic if and 
only if $G$ is finitely generated, soluble, and nilpotent-by-finite.
In particular, if the group of polynomial functions is polycyclic, 
then necessarily it is nilpotent-by-finite. Furthermore, we prove that $G$ itself 
is polycyclic if and only if the subgroup of polynomial 
functions which send $(1,\ldots,1)$ to $1$ is finitely generated and 
soluble.
\end{abstract}
\maketitle
%
%%%%%%%%% 1
\section{Introduction and Main Results}
Let $G$ be a group and let $F_{n}$ be the free group of rank $n>0$
freely generated by $x_{1},\ldots,x_{n}$. More or less explicitely, 
the free product $G[x_{1},\ldots,x_{n}]:=G \ast F_{n}$ 
frequently occurs in group theory, 
for example in the study of equations in groups. Actually, in the 
class of $G$-groups, $G[x_{1},\ldots,x_{n}]$ plays the role of the polynomial 
ring $K[x_{1},\ldots,x_{n}]$ in the class of $K$-algebras (following 
the terminology used in \cite{BMR}, a $G$-group is by definition a group 
containing a designated copy of $G$). If 
$f(x_{1},\ldots,x_{n})\in G[x_{1},\ldots,x_{n}]$ is a ``polynomial'', 
one can define in an obvious way  the associated polynomial function
$(t_{1},\ldots,t_{n})\to f(t_{1},\ldots,t_{n})$, which is a map from 
$G^n$ to $G$. We shall denote by 
$\overline{G}[x_{1},\ldots,x_{n}]$
the set of polynomial functions (understood: in $n$ variables and with 
coefficients in $G$). When $n=1$, we shall write  
$\overline{G}[x]$ instead of  $\overline{G}[x_{1}]$. We define in a 
natural way the product of two polynomial functions by carrying over 
to $\overline{G}[x_{1},\ldots,x_{n}]$ the product of $G[x_{1},\ldots,x_{n}]$.
In other words, the product of $(t_{1},\ldots,t_{n})\to f(t_{1},\ldots,t_{n})$
and $(t_{1},\ldots,t_{n})\to g(t_{1},\ldots,t_{n})$ is equal to the polynomial 
function $$(t_{1},\ldots,t_{n})\to f(t_{1},\ldots,t_{n})g(t_{1},\ldots,t_{n}).$$ 
Equipped with this operation,  
$\overline{G}[x_{1},\ldots,x_{n}]$ is a group and the map changing
$f(x_{1},\ldots,x_{n})\in G[x_{1},\ldots,x_{n}]$ 
into the polynomial function
$(t_{1},\ldots,t_{n})\to f(t_{1},\ldots,t_{n})$ is an epimorphism.

In \cite{BMR}, the authors define notions of zero divisor, ideal,... in a 
group (more precisely in a $G$-group); they obtain then a set of 
results showing a surprising similarity to algebraic geometry. 
In particular, a notion of ``equationally Noetherian group''
is introduced (see \cite{BMR} for a definition) and an analogue 
of the Hilbert's basis theorem is proposed 
as a conjecture \cite[p.42]{BMR}, . The notion of ``equationally Noetherian group'' 
is different from the usual notion of Noetherian group.
Recall here that a  Noetherian group is a group satisfying the maximal 
condition for its subgroups; a polycyclic group is a soluble 
Noetherian group.  
The aim of this paper is to investigate the connections between a group
$G$ and the group $\overline{G}[x_{1},\ldots,x_{n}]$ of its polynomial 
functions for the property of polycyclicity. 
Remark in passing that such a question is not really interesting for the group 
$G[x_{1},\ldots,x_{n}]$. Indeed, if $A$ and $B$ are groups such that 
$A\neq \{ 1\}$ and $|B|\geq 3$, the free product $A\ast B$ 
contains a free subgroup of rank 2 \cite[p.177]{LS} and so is not Noetherian.
Using this result, it is not hard to see that $G[x_{1},\ldots,x_{n}]$ 
is Noetherian if and only if $G$ is trivial and $n=1$.

First we characterize the groups $G$ such that 
$\overline{G}[x_{1},\ldots,x_{n}]$ is polycyclic. Notice that in this 
case, our characterization shows that necessarily $\overline{G}[x_{1},\ldots,x_{n}]$ 
is then nilpotent-by-finite.
\begin{thm} 
For any group $G$ and for any positive integer 
$n$, the following assertions are equivalent:

\noindent {\rm (i)} $\overline{G}[x_{1},\ldots,x_{n}]$ is 
polycyclic;

\noindent {\rm (ii)} $G$ is finitely generated, soluble, and 
nilpotent-by-finite;

\noindent {\rm (iii)} $\overline{G}[x_{1},\ldots,x_{n}]$ is finitely generated, 
soluble, and nilpotent-by-finite.
\end{thm}
An immediate consequence of this theorem is the following.
\begin{cor}
 Let $G$ be a group such that 
$\overline{G}[x]$ is polycyclic; then  so is
$\overline{G}[x_{1},\ldots,x_{n}]$ for any positive integer $n$.
\end{cor}

Theorem 1.1 shows in particular that if $G$ is polycyclic, then 
$\overline{G}[x_{1},\ldots,x_{n}]$ is not necessarily polycyclic.
The next theorem characterizes in terms of polynomial functions the 
case where $G$ is polycyclic. Before to state this result, introduce 
a subset of polynomial functions: we denote by $\overline{G}_{1}[x_{1},\ldots,x_{n}]$ 
the set of polynomial functions which send the $n$-tuple $(1,\ldots,1)$ to $1$. This set is 
obviously a normal subgroup of $\overline{G}[x_{1},\ldots,x_{n}]$.
More precisely, it is easy to see that $\overline{G}[x_{1},\ldots,x_{n}]$
is the (internal) semidirect product of 
$\overline{G}_{1}[x_{1},\ldots,x_{n}]$ and the subgroup of constant 
polynomial functions (isomorphic to $G$). We can now state our result.
\begin{thm}
Let $G$ be a finitely generated group and let 
$n$ be a positive integer. Then $G$ is a polycyclic group if and only if 
$\overline{G}_{1}[x_{1},\ldots,x_{n}]$ is a finitely generated
soluble group.
\end{thm}

Notice that in this theorem, it is necessary to assume that $G$ is 
finitely generated: for example, if $G$ is an abelian group which is 
not finitely generated, then $\overline{G}_{1}[x_{1},\ldots,x_{n}]$ 
is finitely generated and abelian but $G$ is not polycyclic. 

%%%%%%%%%%%%%%%%% 2
\section{Proof of Theorem 1.1}

Consider a group $G$ and an element $v=(a_{1},\ldots,a_{n})\in G^n$.
The map $\phi_{v}: \overline{G}[x_{1},\ldots,x_{n}] \to G$ defined by
$\phi_{v}(f)=f(v)$ (for all $f\in \overline{G}[x_{1},\ldots,x_{n}]$) 
is clearly a homomorphism. Moreover, the intersection 
$\displaystyle \bigcap_{v\in G^n}^{}\ker \phi_{v}$ is trivial. 
It follows that any law of $G$ is also a law for $\overline{G}[x_{1},\ldots,x_{n}]$.
Conversely, since $\overline{G}[x_{1},\ldots,x_{n}]$ contains a copy 
of $G$ (the subgroup of constant polynomial functions), 
each law of $\overline{G}[x_{1},\ldots,x_{n}]$ is a law for $G$.
Thus we can state:

\begin{lem}
For any group $G$ and for any positive integer 
$n$, the variety generated by $G$ coincide with the variety generated by 
$\overline{G}[x_{1},\ldots,x_{n}]$ (that is, $G$ and
$\overline{G}[x_{1},\ldots,x_{n}]$ have the same set of laws).
\end{lem}

Now we introduce some notations used in the statement of the next lemma. 
Let $a,b$ be elements of a group $G$. 
As usual, $[a,\,_{k}b]$ is defined for each integer $k\geq 0$ by 
$[a,\,_{0}b]=a$ and $[a,\,_{k+1}b]=[[a,\,_{k}b],b]$ (where 
$[a,b]=a^{-1}b^{-1}ab$). We shall write $\langle a^{\langle a,b \rangle} 
\rangle$ for the normal closure of $a$ in the subgroup generated by 
$a$ and $b$ and  $\langle a^{\langle a,b \rangle} 
\rangle '$ for its derived subgroup. 
Suppose that there exists a relation of the form
$$w(a,b)[a,\ _{r}b]^{e_0}[a,\ _{r+1}b]^{e_1}\ldots [a,\ _{r+s}b]^{e_{s}}=1$$
with $r,s \in {\mathbb N}$, $e_0,e_1,\ldots ,e_s\in {\mathbb Z}$ ($e_0, 
e_{s}\neq 0$)
and $w(a,b)\in \langle a^{\langle a,b \rangle} \rangle '$. 
We denote by $\Omega_{\star}(a,b)$ (respectively $\Omega^{\star}(a,b)$) the
least integer $\vert e_0\vert$ (respectively $\vert e_s\vert$)
with this property. If $a$ and $b$ do not satisfy a relation of the previous
form, we set $\Omega_{\star}(a,b)=\Omega^{\star}(a,b)=+\infty$. In this 
way, $\Omega_{\star}$ and $\Omega^{\star}$ are two functions 
from $G^2$ to the set ${\mathbb N}^{*}\cup \{ +\infty \}$. 
We proved in \cite{EN} the following result:
 
\begin{lem}
\cite[Corollary 1]{EN} Let $G$ be a finitely generated soluble group.
Then  $G$ is nilpotent-by-finite if and only if for any $a,b\in G$, 
$\Omega^{\star}(a,b)=1$ and the sequence 
$\bigl(\Omega_{\star}(a,b^k)\bigr)_{k>0}$ is bounded.
\end{lem}

\noindent {\em Proof of Theorem 1.1.} (i)$\Rightarrow$(ii). Consider a group $G$ 
such that $\overline{G}[x_{1},\ldots,x_{n}]$ is 
polycyclic for some positive integer $n$. Actually, since plainly
$\overline{G}[x_{1},\ldots,x_{n}]$ contains a copy of
$\overline{G}[x]$, we can assume that $n=1$. First remark that
$G$ is a finitely generated soluble group for $\overline{G}[x]$ contains a 
copy of $G$. It remains to prove that $G$ is 
nilpotent-by-finite. Let $a$ be an element 
in $G$. Consider the subgroup $H\leq \overline{G}[x]$ generated by 
the polynomial functions $f_{k}:t\to [a,\,_{k} t]$, for all positive 
integers $k$.
In fact, there exists  an integer $m>0$ such that $f_{1},\ldots,f_{m}$ 
generates $H$ for $\overline{G}[x]$ is polycyclic. In particular, we 
can write $f_{m+1}$ as a product of factors of the form $f_{j}^{\epsilon_{j}}$
($j=1,\ldots,m$, $\epsilon_{j}\in {\mathbb Z}$). From this writing, we deduce 
a relation of the form
$$f_{m+1}=w(a,x)f_{1}^{e_{1}}\ldots f_{m}^{e_{m}} \;\;\; (e_{j}\in {\mathbb 
Z}),$$ where $w(a,x)$ belongs to the derived subgroup of $H$. It 
follows in $G$ the equality
$$[a,\,_{m+1} b]=
w(a,b)[a,\, b]^{e_{1}}[a,\,_{2} b]^{e_{2}}\ldots [a,\,_{m} b]^{e_{m}} 
\;\;\; ({\rm for \: all}\: b\in G),$$
with $w(a,b)\in \langle a^{\langle a,b \rangle} \rangle '$.
Notice that this relation is independant of $b$.
Hence $\Omega^{\star}(a,b)=1$ and the sequence 
$\bigl(\Omega_{\star}(a,b^k)\bigr)_{k>0}$ is bounded. 
Therefore, by Lemma 2.2, $G$ is nilpotent-by-finite.

\noindent (ii)$\Rightarrow$(iii). Now suppose that the group $G$ is finitely generated, 
soluble and  nilpotent-by-finite. More precisely, suppose that $G$ is 
(nilpotent of class $\nu$)-by-(exponent $\epsilon$), soluble of 
derived length $\rho$, and generated by $g_{1},\ldots ,g_{d}$. Then, 
by Lemma 2.1, $\overline{G}[x_{1},\ldots,x_{n}]$ is also 
(nilpotent of class $\nu$)-by-(exponent $\epsilon$) and soluble of 
derived length $\rho$. Besides, $\overline{G}[x_{1},\ldots,x_{n}]$ is 
finitely generated: it is plain that this group is generated by the 
$d$ constant functions $(t_{1},\ldots,t_{n})\to g_{i}$ 
($i=1,\ldots,d$) and by the $n$ functions $(t_{1},\ldots,t_{n})\to 
t_{j}$ ($j=1,\ldots,n$).

\noindent (iii)$\Rightarrow$(i). This last implication is an immediate 
consequence of well known results.
\null\hfill $\square$ 

%%%%%%%%%%%%%%%%% 3
\section{Proof of Theorem 1.2}

Let $\varphi$ be an automorphism of a finitely 
generated abelian group $A$ (written additively).
Then, by a result of Cohen \cite{CO}, one can find a 
monic polynomial $P\in {\mathbb Z}[T]$ with constant term $1$ such 
that $P(\varphi)=0$ (see also \cite[Theorem 2]{HR}). If 
$P=T^{\lambda}+\epsilon_{\lambda-1}T^{\lambda-1}+
\cdots +\epsilon_{2}T^2 +\epsilon_{1}T+1$, 
we can state this 
result with the multiplicative notation under the form: 

\begin{lem}
\cite{CO} Let $\varphi$ be an automorphism of a finitely 
generated abelian group $A$ (written multiplicatively). Then there exist
a positive integer $\lambda$ and integers 
$\epsilon_{1},\ldots, \epsilon_{\lambda-1}$  such 
that 
$$\varphi^{\lambda}(a)\varphi^{\lambda-1}(a)^{\epsilon_{\lambda-1}}
\ldots \varphi^2(a)^{\epsilon_{2}} \varphi(a)^{\epsilon_{1}}a=1$$
for all $a\in A$.
\end{lem}

The next lemma is an extension of Lemma 3.1 to polycyclic groups.
\begin{lem} 
Let $\varphi$ be an automorphism of a polycyclic group
$G$. Then there exist positive integers $\mu_{1}, \ldots,\mu_{k}$ 
(with $\mu_{i}<\mu_{k}$ for $i=1,\ldots,k-1$) and integers 
$\eta_{1},\ldots, \eta_{k-1}$  such that 
$$\varphi^{\mu_{k}}(t)\varphi^{\mu_{k-1}}(t)^{\eta_{k-1}}
\ldots \varphi^{\mu_{2}}(t)^{\eta_{2}} 
\varphi^{\mu_{1}}(t)^{\eta_{1}}t=1$$
for all $t\in G$.
\end{lem}

\begin{proof}
We argue by induction on the derived length $\rho$ of $G$.
For $\rho =1$, the result is given by Lemma 3.1. Now consider the case
$\rho >1$ and suppose that the result holds for $\rho -1$. 
Put $A=G^{(\rho -1)}$. Since the subgroup $A$ is characteristic,
$\varphi$ induces in $G/A$ an automorphism. By the 
inductive hypothesis applied to $G/A$ with this automorphism,
there exist positive integers $\mu_{1}, \ldots,\mu_{k}$ 
($\mu_{i}<\mu_{k}$ for $i=1,\ldots,k-1$) and integers 
$\eta_{1},\ldots, \eta_{k-1}$  such that 
$$\varphi^{\mu_{k}}(t)\varphi^{\mu_{k-1}}(t)^{\eta_{k-1}}
\ldots \varphi^{\mu_{2}}(t)^{\eta_{2}} 
\varphi^{\mu_{1}}(t)^{\eta_{1}}t$$
belongs to $A$ for all $t\in G$. Now notice that $\varphi$ defines by 
restriction an automorphism in $A$. Apply Lemma 3.1 
to $A$ with this automorphism; thus there exist
a positive integer $\lambda$ and integers 
$\epsilon_{1},\ldots, \epsilon_{\lambda-1}$  such 
that 
$$\varphi^{\lambda}(a)\varphi^{\lambda-1}(a)^{\epsilon_{\lambda-1}}
\ldots \varphi^2(a)^{\epsilon_{2}} \varphi(a)^{\epsilon_{1}}a=1$$
for all $a\in A$. Clearly, by taking
$$a=\varphi^{\mu_{k}}(t)\varphi^{\mu_{k-1}}(t)^{\eta_{k-1}}
\ldots \varphi^{\mu_{2}}(t)^{\eta_{2}} 
\varphi^{\mu_{1}}(t)^{\eta_{1}}t$$ in this relation (for any $t\in G$), 
we obtain the required result. 
\end{proof}
   
\begin{lem}  Let $\varphi$ be an automorphism of a polycyclic group
$G$. Then there exists a positive integer $\xi$ such that,
for each integer $\zeta\geq \xi$ (resp. $\zeta\leq -\xi$),
there exist a positive integer $m$ and integers 
$\xi_{1}, \ldots,\xi_{m},\theta_{1},\ldots, \theta_{m}$ with 
$0\leq \xi_{i}<\xi$ (resp. $-\xi < \xi_{i}\leq 0$) for each $i=1,\ldots,m$,
such that 
$\varphi^{\zeta}(t)=
\varphi^{\xi_{1}}(t)^{\theta_{1}}
\ldots \varphi^{\xi_{m}}(t)^{\theta_{m}} $
for all $t\in G$.
\end{lem} 

\begin{proof}
Set $\xi=\mu_{k}$, where $\mu_{k}$ is the integer defined 
in Lemma 3.2. If $\zeta=\xi$, Lemma 3.2 gives the required relation, namely
$$\varphi^{\xi}(t)=t^{-1}\varphi^{\mu_{1}}(t)^{-\eta_{1}}\ldots
\varphi^{\mu_{k-1}}(t)^{-\eta_{k-1}}\;\;\; ({\rm for\: all}\: t\in G).$$ By 
using this last relation, an easy induction prove the property for all 
$\zeta\geq\xi$. When  $\zeta\leq -\xi$, the argument is similar, but one 
use the relation 
$$ \varphi^{-\xi}(t)=\varphi^{-\mu_{k}}(t)=\varphi^{\mu_{1}-\xi}(t)^{-\eta_{1}}
\varphi^{\mu_{2}-\xi}(t)^{-\eta_{2}} \ldots 
\varphi^{\mu_{k-1}-\xi}(t)^{-\eta_{k-1}}t^{-1},$$ 
which follows from Lemma 3.2.
\end{proof}
 
The proof of the next lemma will be omitted: 
the first part is a direct consequence of Lemma 3.3 and 
the second part follows from the first part by induction on $d$. 
 
\begin{lem} 
Let $G$ be a polycyclic group. Then:

\noindent {\rm (i)} If $\varphi$ is an automorphism of $G$, 
there exists a positive integer $\xi$ such that,
for each $\zeta\in {\mathbb Z}$, there exist a positive integer $m$ and integers 
$\xi_{1}, \ldots,\xi_{m},\theta_{1},\ldots, \theta_{m}$, 
with $|\xi_{i}| <\xi$ for $i=1,\ldots,m$, 
such that 
$\varphi^{\zeta}(t)=
\varphi^{\xi_{1}}(t)^{\theta_{1}}
\ldots \varphi^{\xi_{m}}(t)^{\theta_{m}} $
for all $t\in G$;

\noindent {\rm (ii)} More generally, if $\varphi_{1},\ldots,\varphi_{d}$ are $d$ 
automorphisms of $G$, 
there exists a positive integer $\xi$ such that,
for each $(\zeta_{1},\ldots,\zeta_{d})\in {\mathbb Z}^d$, 
there exist a positive integer $m$ and integers 
$\xi_{1,j}, \ldots,\xi_{m,j}$ ($j=1,\ldots,d$), $\theta_{1},\ldots, \theta_{m}$, 
with $|\xi_{i,j}| <\xi$ for $i=1,\ldots,m$ and $j=1,\ldots,d$,
such that 
$$\varphi_{1}^{\zeta_{1}}\circ \ldots \circ \varphi_{d}^{\zeta_{d}}(t)=
\left( \varphi_{1}^{\xi_{1,1}}\circ \ldots \circ 
\varphi_{d}^{\xi_{1,d}}(t)\right)^{\theta_{1}}
\ldots 
\left( \varphi_{1}^{\xi_{m,1}}\circ \ldots \circ 
\varphi_{d}^{\xi_{m,d}}(t)\right)^{\theta_{m}}$$
for all $x\in G$.
\end{lem} 

Finally, as a particular case, we obtain the following result, 
essential in the proof of Theorem 1.2. 

\begin{lem}
Let $b_{1},\ldots,b_{d}$ be $d$ fixed elements of a 
polycyclic group $G$ and let $H$ be the subgroup of $\overline{G}[x]$ 
generated by the 
polynomial functions (in one variable) of the form
$t\to b_{d}^{-\beta_{d}} \ldots b_{1}^{-\beta_{1}} t 
b_{1}^{\beta_{1}}\ldots b_{d}^{\beta_{d}}$, 
with $\beta_{1},\ldots,\beta_{d}\in {\mathbb Z}$. 
Then there exists a positive integer $\xi$ such that $H$ is generated 
by the polynomial functions of the form
$$t\to b_{d}^{-\beta_{d}} \ldots b_{1}^{-\beta_{1}} t 
b_{1}^{\beta_{1}}\ldots b_{d}^{\beta_{d}},$$
where $\beta_{1},\ldots,\beta_{d}$ are integers such that
$|\beta_{i}| <\xi$ (in other words, $H$ is finitely generated).
\end{lem}

\begin{proof}
It suffices to apply Lemma 3.4(ii) in the case where  
$\varphi_{j}$ is the inner automorphism of $G$ defined 
by $\varphi _{j} (t)=b_{j}^{-1}tb_{j}$, with $j=1,\ldots,d$.
\end{proof}

We need again two easy lemmas:

\begin{lem} 
In each polycyclic group, there exists a finite 
sequence $b_{1},\ldots,b_{d}\in G$ such that any element $a\in G$
may be written in the form 
$a=b_{1}^{\beta_{1}}\ldots b_{d}^{\beta_{d}}$, where
$\beta_{1},\ldots,\beta_{d}$ are integers.
\end{lem}

\begin{proof}
The group $G$ being polycyclic, it has a series
$$\{ 1\} =G_{d+1}\unlhd G_{d}\unlhd\cdots\unlhd G_{2}\unlhd G_{1}=G$$
in which each factor $G_{j}/G_{j+1}$ is cyclic. For $j=1,\ldots,d$, 
choose in each $G_{j}$ an element  $b_{j}$ such that $b_{j}G_{j+1}$
generates  $G_{j}/G_{j+1}$. Then clearly the sequence 
$b_{1},\ldots,b_{d}$ satisfies the statement of the lemma.
\end{proof}

\begin{lem} 
Let $G$ be a finitely generated group in 
which the normal closure of each element is finitely generated and soluble.
Then $G$ is polycyclic.
\end{lem}

\begin{proof}
Clearly, a finitely generated group 
$\langle g_{1},\ldots,g_{d}\rangle$ such that 
the normal closure of each generator $g_{j}$ is soluble is 
soluble itself. Thus $G$ is soluble.
Now suppose that $G$ is not polycyclic. 
Since a polycyclic group is finitely presented, we can assume that 
each proper homomorphic image of $G$ is polycyclic (see for example 
\cite[Part 2, Lemma 6.17]{RO}). The group $G$ being soluble, it contains 
a non-trivial abelian normal subgroup; let $A$ be the normal closure 
of a non-trivial element of this subgroup. Then $G/A$ is polycyclic.
Furthermore, $A$ is abelian and finitely generated by hypothesis; 
thus $A$ is polycyclic. It follows that $G$ is polycyclic, a 
contradiction. 
\end{proof}

\noindent {\em Proof of Theorem 1.2.} First suppose that $G$ is a polycyclic group.
Since $G$ is soluble, Lemma 2.1 shows that $\overline{G}[x_{1},\ldots,x_{n}]$
is soluble too; thus so is the subgroup $\overline{G}_{1}[x_{1},\ldots,x_{n}]$.
For the second part of the property,
we begin with the case $n=1$; thus we want prove that 
$\overline{G}_{1}[x]$ is finitely generated. It is not difficult to see 
that this group is generated by the functions of the form
$t\to a^{-1}ta$, with $a\in G$. By Lemma 3.6, there exist elements
$b_{1},\ldots,b_{d}\in G$ depending only on $G$ such that each $a\in G$
may be written in the form $a=b_{1}^{\beta_{1}}\ldots b_{d}^{\beta_{d}}$.
Therefore $\overline{G}_{1}[x]$ is generated by the functions of the form
$t\to b_{d}^{-\beta_{d}}\ldots b_{1}^{-\beta_{1}} t
b_{1}^{\beta_{1}}\ldots b_{d}^{\beta_{d}}$ ($\beta_{j}\in {\mathbb Z}$).
We can deduce then from Lemma 3.5 that the group $\overline{G}_{1}[x]$ is finitely generated.

\noindent Now suppose that $n$ is an arbitrary positive integer.
For each $i\in\{ 1,\ldots,n\}$, one can define a monomorphism
$\Psi_{i}:\overline{G}_{1}[x]\to 
\overline{G}_{1}[x_{1},\ldots,x_{n}]$ in the following way: if 
$f\in \overline{G}_{1}[x]$, then 
$\Psi_{i} (f) (t_{1},\ldots,t_{n})=f(t_{i})$ for all 
$(t_{1},\ldots,t_{n})\in G^n$. Therefore each subgroup
$\Psi_{i}\left( \overline{G}_{1}[x] \right)$ is isomorphic to
$\overline{G}_{1}[x]$ and so is finitely generated.
Now remark that $\overline{G}_{1}[x_{1},\ldots,x_{n}]$ 
is generated by the functions of the form
$(t_{1},\ldots,t_{n})\to a^{-1}t_{i}a$ (with $i=1,\ldots,n$ and $a\in G$);
furthermore, such a function belongs to 
$\Psi_{i}\left( \overline{G}_{1}[x] \right)$. Thus 
$\overline{G}_{1}[x_{1},\ldots,x_{n}]$ is generated by 
$\Psi_{1}\left( \overline{G}_{1}[x] \right)\cup \ldots \cup
\Psi_{n}\left( \overline{G}_{1}[x] \right)$ and so is finitely 
generated, as required.

\noindent Conversely, suppose now that 
$\overline{G}_{1}[x_{1},\ldots,x_{n}]$ is a finitely generated
soluble group.
Consider an element $v=(a,1,\ldots,1)\in G^n$, where $a$ is an 
element of $G$.
The map $\Phi_{v}: \overline{G}_{1}[x_{1},\ldots,x_{n}] \to G$ defined by
$\Phi_{v}(f)=f(v)$ is clearly a homomorphism. Moreover, it is easy to 
see that  $\Phi_{v} \left( \overline{G}_{1}[x_{1},\ldots,x_{n}]  
\right)$ coincide with the normal closure of $a$ in $G$. Thus the 
normal closure of each element in $G$ is a finitely generated soluble 
subgroup. Since $G$ is finitely generated by hypothesis, we can apply Lemma 
3.7, and hence $G$ is polycyclic. This completes the proof of the theorem.
\null\hfill $\square$ 

%
%%%%%% References
%

%

\end{document}